\documentclass[amscd]{amsart}

\newtheorem{Thm}{Theorem}

\newtheorem{Lem}{Lemma}
\newtheorem{Prop}{Proposition}

\newtheorem{Claim}{Claim}
\newtheorem{Conj}{Conjecture}

\theoremstyle{remark}
\newtheorem{Rem}{Remark}
\newtheorem{Def}{Definition}
\newtheorem{Ex}{Example}

\newcommand{\F}{{\mathcal F}}
\newcommand{\A}{{\mathcal A}}
\renewcommand{\O}{{\mathcal Z}}
\newcommand{\MM}{{\mathcal M}}

\newcommand{\Gm}{{{\mathbb G}_m}}

\renewcommand{\AA}{{\mathcal A}}

\newcommand{\g}{{\mathfrak g}}

\newcommand{\h}{{\mathfrak h}}

\newcommand{\BB}{{\mathcal B}}
\newcommand{\CC}{{\mathcal C}}
\newcommand{\FF}{{\mathcal F}}
\newcommand{\GG}{{\mathcal G}}
\newcommand{\OO}{{\mathcal O}}
\newcommand{\PP}{{\mathcal P}}

\newcommand{\Zet}{{\mathbb Z}}

\newcommand{\imbed}{\hookrightarrow}

\newcommand{\RE}{{\mathbb R}}
\newcommand{\Ce}{{\mathbb C}}

\newcommand{\<}{\langle}
\renewcommand{\>}{\rangle}

\newcommand{\Alc}{{{\mathrm{Alc}} }}

\newcommand{\bla}{{\boldsymbol{\lambda}}}

\def\square{\hbox{\vrule\vbox{\hrule\phantom{o}\hrule}\vrule}}

\newcommand{\al}{\alpha}
\newcommand{\la}{\lambda}
\newcommand{\La}{\Lambda}

\newcommand{\ii}{\sqrt{-1}}

\newcommand{\wt}{\widetilde}

\newcommand{\epf}{\square}

\begin{document}

\title[(Real) stabilities conditions for slices]{Stability conditions for Slodowy slices and real variations of stability}
%\author{Rina Anno, Roman Bezrukavnikov and Ivan Mirkovi\' c}

\author{    Rina Anno     }
\address{\small
Department of Mathematics,
University of Pittsburg, Pittsburgh, PA 15260, USA
}
\email{                anno@pitt.edu        }

\author{
Roman Bezrukavnikov
}
\address{\small
Department of Mathematics, Massachusetts Institute of Technology, 77 Massachusetts ave.,
Cambridge, MA 02139, USA
}
\address{\small
National Research University Higher School of Economics,
International Laboratory of Representation
Theory and Mathematical Physics,
20 Myasnitskaya st., Moscow 101000, Russia}
\email{
bezrukav@math.mit.edu
}
\author{    Ivan Mirkovi\'c     }
\address{\small
Department of Mathematics and Statistics,
University of Massachusetts, Amherst, MA 01003, USA
}
\email{                mirkovic@math.umass.edu        }

\begin{abstract} The paper provides new examples of an explicit submanifold in
 Bridgeland stabilities space of a local Calabi-Yau.

More precisely, let $X$ be the standard resolution of a transversal slice to an adjoint nilpotent orbit of a simple Lie algebra over $\Ce$. An action of the affine braid group on the derived category $D^b(Coh(X))$ and a collection of $t$-structures on this category permuted by the action have been constructed in \cite{BR} and \cite{BM} respectively.
 In this note we show that the $t$-structures come from points in a certain connected submanifold in the space of Bridgeland stability conditions.
 The submanifold is a covering of a submanifold in the dual space to the Grothendieck group, and the affine braid group acts by
 deck transformations.

 We also propose a new variant of definition of stabilities on a triangulated category, which we call a "real variation of stability conditions" and discuss its relation to Bridgeland's definition. The main theorem provides an
illustration of such a relation. We state a conjecture by the second author and A.~Okounkov
on examples of this structure arising from symplectic resolutions of singularities and its relation
to equivariant quantum cohomology.
We  verify this conjecture in our examples.

\end{abstract}

\maketitle

\centerline{\em{To Borya Feigin with gratitude and best wishes 
on his anniversary}}
\section{Introduction}
The goal of this work is twofold. One aim, achieved in Theorem \ref{thm1} is to describe new examples of an explicit connected submanifold in the space of locally finite Bridgeland stabilities on the derived category of a local
Calabi-Yau manifold $X$. The examples in question come from a simple algebraic group $G$,
more precisely $X$ is the resolution of a transversal slice to an adjoint orbit in the nilpotent cone
of a simple algebraic group. It is well know that minimal resolutions of Kleinian surface singularities 
belong to this class, in that special case we recover a  weaker variant of  Bridgeland's result
 \cite{Br}, see Example \ref{ExKl}. 
Notice that the real dimension of the submanifold in the space of stabilities
 we describe 
%in Theorem \ref{thm1}  
is twice the second Betti number
of the flag variety of $G$; in almost all cases (see footnote to Theorem \ref{thm2} for more precise information)
 this is equal to twice the second Betti number of $X$. It is our understanding that mirror symmetry conjectures suggest existence of a canonical submanifold in the stability space $Stab(D^b(Coh(M)))$ for a Calabi-Yau manifold $M$, called the stringy moduli space, whose
 real dimension is $2b_2(M)$: under mirror duality it should correspond to a covering of the moduli space of deformations of the dual Calabi-Yau manifold. 
It seems natural to expect that the submanifold we describe  is related to the stringy moduli space of $X$.

The second goal is to explore the parallelisms between  structures discovered, respectively, by representation theorists working in Kazhdan-Lusztig theory and its generalizations and by algebraic geometers studying   Bridgeland stability conditions. The first point to mention here is the
action of the braid group (or its generalizations) on the derived category of modules, which
is the key ingredient in Kazhdan-Lusztig theory. The braid group can be described as the 
fundamental group of a space which, in our opinion, should be thought of as the counterpart of the space $\mathcal K$ of Kaehler parameters in algebraic geometry, the action of the fundamental
group  $\pi_1({\mathcal K})$ on the derived category of coherent sheaves has been constructed
in a number of examples (see references in section \ref{secrv}). Localization theory in positive
characteristic \cite{BMR} relates derived category of modules to that of coherent sheaves,
in this setting the two constructions can be directly related as explained below.

Furthermore, as the proof of the Theorem \ref{thm1} shows, in some examples a variant of the 
{\em central
charge} map appearing in Bridgeland's definition arises naturally from polynomials describing {\em dimensions of modules} in positive characteristic.
These examples led us to
introduce new definitions of structures related to but different from
Bridgeland stability manifold: these are   (symmetric) real variations of stability conditions
and  local systems of categories with stabilities. 

We discuss heuristic relation of our definitions
to Bridgeland's one in Remark \ref{discu}, while comparing Theorem \ref{thm1} 
to Theorem \ref{thm2} one gets an example of a
precise connection between the two. We also state a conjecture (Conjecture \ref{conj}
in section \ref{secsymp}) due to  A.~Okounkov and the second author on existence of such a structure on derived categories of general symplectic resolutions
and their relation to equivariant quantum cohomology. We explain in Theorem \ref{thm2} 
that validity of this conjecture in our examples follows from known results.

\subsection{Acknowledgements} We thank Andrei Okounkov for the permission to
quote a conjecture partly due to him and 
Tom Bridgeland for help with the proof of Lemma \ref{LBr}. The second author thanks Leonid Rybnikov for discussions
leading to conjecture in Remark \ref{new_conj}.
 A part of the work on the paper was done while the last two authors enjoyed the hospitality of Institute for Advanced Study at Hebrew University of Jerusalem,
they would like to thank that Institution for excellent work conditions. 
R.B. and I.M. were supported by NSF grants.

We are grateful for the opportunity to dedicate this work to 
Borya Feigin; his unique  vision of things in and beyond  mathematics
%, of raising mathematicians
 % and of other related matters 
has been an important influence and a source of inspiration for the authors
at various stages of their professional life. 

\section{Bridgeland stabilities for the resolutions of transversal slices}\label{sec1}
Let $\g$ be a %semi-
simple Lie algebra over $\Ce$, let $e\in \g$ be a non-principal nilpotent element,
$Y\subset \g$ be a transversal (Slodowy) slice to the $G$-orbit of $e$. Let $\BB=G/B$ be the flag variety of $G$, and $\pi:T^*(\BB)\to \g$ be the Springer (moment) map. Set $\BB_e=\pi^{-1}(e)$ and  $X=\pi^{-1}(Y)$. Set $\CC=D^b(Coh_{\BB_e}(X))$ where
$Coh_{\BB_e}(X)$ is the category of coherent sheaves on $X$ supported on $\BB_e$.

Certain $t$-structures on $\CC$ were constructed in \cite{BM} (announced in \cite{ICM}). In this note we show that they arise from an explicit connected subset of the space $Stab(\CC)$ of locally finite Bridgeland stability conditions on $\CC$. To state the result we need more notations.

Let $\h$ denote the (abstract) Cartan algebra of $\g$.

We have $\h^*=\Lambda\otimes _\Zet \Ce$ where $\Lambda$ is the weight lattice. Let
$\h^*_\RE=\Lambda \otimes_\Zet \RE\subset \h^*$ be the real dual Cartan. The affine Weyl group $W_{aff}$ acts on $\h^*$ and on $\h^*_\RE$
by affine-linear transformations.
Let $\h^*_{reg}$ be the union of free orbits of $W_{aff}$ on $\h^*$, thus $\h^*_{reg}$ is the
 complement to the affine coroot hyperplanes $H_{\check{\al}, n}=\{\la\in \h^*\ | \ \langle \la , \check{\al} \rangle =n\}$, where $n\in \Zet$ and $\check{\al}$ is a coroot.

\subsection{Action of the affine braid group and $t$-structures assigned to alcoves.} \label{11}
Recall that an {\em alcove} is a connected component of $\h^*_{\RE,reg}=\h^*_{reg}\cap \h^*_\RE$. The natural action of the affine Weyl group on $\h^*$ induces a simply transitive action of $W_{aff}$ on the set of alcoves. We denote this set by $Alc$.

The argument below is based on the construction of \cite{BM} which assigns a $t$-structure $\tau_A$ on $\CC$ to an alcove $A\in Alc$. The $t$-structure $\tau_A$ can be described using
derived localization over a field of characteristic $p>0$ \cite{BMR}. Roughly speaking, modules for
the sheaf $D_\lambda(\BB)$ of twisted differential operators on $\BB$ are closely related to
coherent sheaves on $T^*(\BB)$; on the other hand, the derived category $D^b(D_\lambda(\BB))$
is identified with the derived category of an appropriate quotient of the enveloping algebra
$U\g$. Thus one can get a $t$-structure on $D^b(Coh(T^*(\BB)))$ which is compatible with the tautological $t$-structure on $D^b(U\g-mod)$. The $t$-structure $\tau_A$ arises this way
when the twisting parameter $\lambda$ satisfies the condition $\frac{\lambda+\rho}{p}\in A$.
There exists also a more direct construction of the $t$-structure $\tau_A$ over a characteristic
zero field,  though available proof  of its properties relies on positive characteristic picture.

For future reference we mention that $\tau_A$ is compatible with a $t$-structure
on $D^b(Coh(X))$ which corresponds to the tautological $t$-structure under an equivalence
$D^b(Coh(X))\cong D^b(R-mod_{fg})$ where $R$ is a certain $\O(Y)$-ring which is finite 
as an $\O(Y)$ module
and $R-mod_{fg}$ denotes the category of finitely generated $R$-modules. It follows that the heart of $\tau_A$ is a finite length abelian category.

Let $B_{aff}=\pi_1(\h^*_{reg}/W_{aff})$ be the affine braid group (this is the affine braid group of Langlands dual group in the standard terminology). An action of
$B_{aff}$ on $\CC$ was defined in \cite{BR}. This action permutes the $t$-structures
$\tau_A$. More precisely, to each pair of alcoves $A,\, A'\in Alc$ one can assign an element $b_{A,A'}\in B_{aff}$; it is then shown in \cite{BM} that $b_{A,A'}$ sends
$\tau_A$ to $\tau_{A'}$. To define $b_{A,A'}$ notice that  an element in $B_{aff}$
is determined by a homotopy class of a path connecting two alcoves in $\h^*_{reg}$.
The element $b_{A,A'}$ corresponds to a path $\phi: [0,1]\to \h^*_{reg}$ such that
$\phi(0)\in A$, $\phi(1)\in A'$ and $\phi(t)\in \h^*_\RE + i (\h^*_\RE)^+ $ for $t\in (0,1)$; here $(\h^*_\RE)^+\subset \h^*_\RE$ is the dominant Weyl chamber. This requirement characterizes the homotopy class of $\phi$ uniquely.

For future reference we fix a universal covering $\wt{\h^*_{reg} }$. We also fix a continuous lifting of each alcove $A\in Alc$ to a subset $\wt{A}$ in $\wt{\h^*_{reg} }$, so that
for each two alcoves $A$, $A'$ a path representing $b_{A,A'}$ lifts to a continuous path connecting $\wt{A}$ to $\wt{A'}$.

\subsection{Embedding $\h^*\to K^0(\CC)^*$ and the "quasi-exponential" map}

We identify $H^*(G/B,\Ce)$ with $K^0(Coh(G/B))\otimes \Ce$ by means of the Chern character map. Notice that the class of the line bundle $\OO(\la)$ attached to $\la\in \Lambda$ corresponds to $\exp(\la)\in H^*(G/B)$ where $\la\in \La$ is considered as an element
in $\h^*=H^2(G/B)$; it is a nilpotent element in the commutative algebra $H^*(G/B)$, so its exponent is well defined.

We have a bilinear pairing $K^0(G/B)\times K^0(\CC)\to \Zet$
given by $([\FF], [\GG])= \chi(pr^*(\FF)\otimes \GG)$. Here $\chi$ stands for Euler characteristic and $pr$ for the projection $T^*(G/B)\to G/B$. This gives a map
$H^*(G/B) \to (K^0(\CC)\otimes \Ce)^*$. We will omit complexification from notation where it is not likely to lead to a confusion, and identify an element in $H^*(G/B)$
with its image in $K^0(\CC)^*$.

We extend the map $\La\to H^*(G/B)\cong K^0(Coh(G/B))\otimes \Ce$, $\la\mapsto \exp(\la)$ to
$\h^*$ as follows. Define the "quasi-exponential" map $E:\h^*=\h^*_\RE\times (\ii\h^*_\RE)\to H^*(G/B)$ by:
%
%$$E:x+\ii y\mapsto  exp(x)(1+\ii [exp(y)-1])= exp(x)(1+\ii (y+y^2/2+\dots)).$$
$$E:x+\ii y\mapsto  \exp(x)(1+\ii \exp(y)).$$ %= exp(x)(1+\ii (y+y^2/2+\dots)).$$
%Notice that the real part of $E(z)$ depends only on the real part of $z$, and same for imaginary parts.

In fact, the map 
\begin{equation}\label{xiy}
x+\ii y\mapsto (x,x+y)
\end{equation}
 is a $W_{aff}$ equivariant isomorphism between $\h^*$ and $(\h^*_\RE)^2$, where $W_{aff}$ acts on $(\h^*_\RE)^2$ diagonally.
Written as a map from $(\h^*_\RE)^2$, the map $E$ takes the form $(\la,\mu)\mapsto \exp(\la) +\ii (\exp(\mu))$.

\begin{Rem}
 A variation of the argument below also works for the map $E(z)=\exp(z)$ (with a less explicit and not necessarily open, though still connected neighborhood $V$ of $(\h^*_\RE)^{ar}$). The proof of the statement involving the above quasi-exponential map
is a bit shorter, so we opted for presenting that version of the result.
\end{Rem}

\begin{Lem}
 The map $E$ is compatible with the $W_{aff}$ action where the action on the source is the standard affine linear action
on $\h^*$, and the one on the target is induced by the $B_{aff}$ action on $D^b(Coh(T^*(\BB))$ from \cite{BM}.
\end{Lem}

 {\em Proof.}  Translations
act on the target by twisting with a line bundle and on the source by shifting the real part, thus
it is easy to deduce that the map is compatible with the action of the lattice of translations. Compatibility with the action of the finite Weyl group $W$ follows from
\cite[Theorem 1.3.2]{BM}. \epf

\subsection{The main result}
Define the "almost regular" part $(\h^*_\RE)^{ar}$ of the real Cartan  as the set of points in $\h^*_\RE$
whose stabilizer in $W_{aff}$ has at most two elements.
For $\la,\mu \in \h^*_\RE$ we will write $\la\preceq \mu$ if $\la$ lies in the closure of
the face which contains $\mu$. Here by a face we mean a stratum of the stratification
of $\h^*_\RE$ cut out by the coroot hyperplanes  (thus alcoves are faces of maximal dimension).
Define a neighborhood $V$ of the set of almost regular
real points  $(\h^*_\RE)^{ar}$ in $\h^*=\h^*_\RE\otimes_\RE \Ce$ by:

 $$V=\{(\la, \mu)\in (\h^*_\RE)^{ar}\times \h^*_\RE \subset
 \h^*_\RE\times \h^*_\RE \cong \h^* \ |\ \la\preceq \mu \bigvee \la \in (\h^*_\RE)^{reg} \},$$
 %(\la \in (\h^*_\RE)^{reg}) \bigvee (\la\in \bar{A}, \mu \in A$ for some $A\in Alc) \}$.
 %
where we used the isomorphism inverse to \eqref{xiy}.

 Thus $V$ is an open neighborhood of $(\h^*_\RE)^{ar}$ in $\h^*$. Let $V^{reg}=V\cap \h^*_{reg}$; we have:

  $V^{reg}=\{(\la, \mu)\in (\h^*_\RE)^{ar}\times \h^*_\RE\ |(\la \in (\h^*_\RE)^{reg}) \bigvee (\la\in \bar{A}, \mu \in A$ for some $A\in Alc) \}$.

Let $\wt{V^{reg}}$ be the preimage of $V^{reg}$ in $\wt{\h^*_{reg}}$.

\begin{Thm}\label{thm1}
There exists a unique map $\iota$ from
$\wt{V^{reg}}$ to the space $Stab(\CC)$ of locally finite Bridgeland stability conditions on $\CC$ such that:

\begin{enumerate}
\item The composed map $Z\circ \iota$, where $Z$ is the projection $Stab\to K^0(\CC)^*$,  coincides with the map $\ii E\circ \pi$ where $\pi$ is the projection $\wt{\h^*_{reg}}\to \h^*_{reg}$.

\item For some (equivalently, for any) $A\in Alc$ and $z\in \wt{A}$ the underlying $t$-structure of the stability $\iota(z)$  coincides with $\tau_A$.
\end{enumerate}

 The map $\iota$ is compatible with the action  of $B_{aff}$ which acts on the source by deck transformations, while the action on the target comes from the action on the category $\CC$.
\end{Thm}

{\em Proof.} Uniqueness of a map $\iota$ satisfying (1) and (2) for some fixed alcove $A$ and
$z\in \wt{A}$ follows from  \cite[Theorem 1.2]{Br1} which asserts that the map $Z$ is a local homeomorphism. It remains to show existence of a $B_{aff}$ equivariant map $\iota$ which satisfies (1) and (2) for all $A$, $z\in \wt{A}$. This will be done in section \ref{sect3}.

\begin{Ex}\label{ExKl} If $\dim (X)=2$, i.e. $e$ is sub-regular, $X$ is well known to be the minimal resolution of a Kleinian singularity. In this case a component of the space
$Stab(\CC)$ was described in \cite{Br} (in the special
case of type $A$ it is shown in \cite{IUU} that this component equals the whole space $Stab(\CC)$). 
It is easy to see that our submanifold is contained in the one described in {\em loc. cit.}
\end{Ex}

\begin{Rem}
It is easy to show that $\pi_1(V^{reg}/W_{aff})$ is a free group with $rank(G)$ generators. This group surjects to $B_{aff}=\pi_1(\h^*_{reg}/W_{aff})$. (The same remains true if $V$ is replaced by any sufficiently small convex $W_{aff}$ invariant neighborhood of $(\h^*_\RE)^{ar}$ in $\h^*$). Thus the covering  $\wt{V^{reg}}\to V^{reg}$ is connected but it is far from being universal. It would be interesting to construct an explicit subset in $Stab(\CC)$ which is a universal covering of a domain whose fundamental group is isomorphic to the affine braid group.
% The difficulty seems to come from the fact that
%for an irreducible object $L$ in the heart of $\tau_A$ the corresponding functional $d_L$ can vanish on several faces of $A$
%(see below).
\end{Rem}

\begin{Rem}
By a standard argument (cf. \cite{Br}, \cite{BrK3}) injectivity of the map $\iota$ is equivalent to 
the fact that the orbit of $\tau_A$ under the action of  $B_{aff}$ on the set of $t$-structures on $\CC$
is free. Notice that the result of \cite{ST} in type $A$ and of \cite{BT} in general 
implies (at least if $G$ is  simply-laced)  that the orbit of $\tau_A$ under the action of the  
subgroup $B\subset B_{aff}$ is free; here $B$ is the Artin braid group associated to the finite
Dynkin graph.
\end{Rem}

\section{Real variations of stabilities}
In this section we discuss the motivation for the main result and suggest some new definitions.
We also describe a conjecture by A.~Okounkov and the second author and verify it in our examples.

\subsection{Real variations of stability conditions: definition}
We expect that the following structure is relevant, at least in some examples, for understanding the  aspects of Calabi-Yau categories which have been studied in the literature via the concept of Bridgeland stabilities.

Let $\CC$ be a finite type triangulated category and $V$ a real vector space.
Suppose that a discrete collection $\Sigma$ of affine hyperplanes in $V$ is fixed, let $V^0$ denote their complement. For each hyperplane in $\Sigma$ consider the parallel hyperplane passing through
zero,
let $\Sigma_{lin}$ be the set of those linear hyperplanes. Fix a component $V^+$ of the complement to the union of hyperplanes in $\Sigma_{lin}$. The choice of $V^+$ determines for each $H\in \Sigma$ the choice of the positive half-space
$(V\setminus H)^+\subset V\setminus H$, where $(V\setminus H)^+=H+V^+$.
By an {\em alcove} we mean a connected component of the complement to hyperplanes in
$\Sigma$ and we let $\Alc$ denote the set of alcoves.
For two alcoves
 $A$, $A'\in \Alc$ sharing a codimension one face which is contained
in a hyperplane $H\in \Sigma$ we will say that $A'$ is {\em above}
$A$ and $A$ is below $A'$ if $A'\in (V\setminus H)^+$.

\begin{Def}\label{def1}
A {\em real variation of stability conditions} on $\CC$ parametrized by $V^0$ and directed to $V^+$
is the data $(Z,\tau)$, where $Z$ (the central charge) is a polynomial map $Z:V\to (K^0(\CC)\otimes \RE)^*$, and $\tau$ is a map from $\Alc$ to the set of bounded
$t$-structures on $\CC$, subject to the following conditions.

\begin{enumerate}
\item
If $M$ is a nonzero object in the heart of $\tau(A)$, $A\in \Alc$, then
$\< Z(x), [M]\> >0$ for $x\in A$.

\item
 Suppose $A$, $A'\in \Alc$ share a codimension one face $H$ and $A'$ is above $A$.
Let $\AA$ be the heart of $\tau(A)$; for $n\in {\mathbb{N}}$
  let $\AA_n\subset \AA$
 be the full
subcategory in $\AA$ given by: $M\in \AA_n$  if the polynomial function on $V$,
$x\mapsto \< Z(x), [M] \>$ has zero of order at least $n$ on $H$.
One can check that $\AA_n$ is a Serre subcategory in $\AA$, thus $\CC_n=\{
C\in \CC\ |\ H^i_{\tau(A)}(C)\in \AA_n\}$ is a thick subcategory in $\CC$.
We require that

\begin{enumerate}
\item The $t$-structure $\tau(A')$ is compatible with the filtration by thick
subcategories $\CC_n$.

\item The functor of shift by $n$ sends the $t$-structure on $gr_n(\CC)=\CC_n/\CC_{n+1}$ induced
by $\tau(A)$ to that induced by $\tau(A')$. In other words, $$gr_n(\AA')=gr_n(\AA)[n]$$ where $\AA'$ is the heart of $\tau(A')$, $gr_n=\AA'_n/\AA'_{n+1}$, $\AA'_n=\AA'\cap \CC_n$.

\end{enumerate}
\end{enumerate}

\end{Def}

\begin{Rem}\label{discu}
 Requirement (1) of the Definition implies that $(\ii Z,\tau)$ define a map from $V^0$ to the space of Bridgeland stabilities $Stab(\CC)$.
  %(recall that the positive ray is contained in the upper half-plane according to Bridgeland's conventions);
Since $V^0$ is disconnected, this structure by itself does not provide any relation between
the different $t$-structures, thus it is too weak to yield interesting results. Axiom (2) connects
 the $t$-structures assigned to different connected components of $V^0$; it is based
 on the same intuition as Bridgeland's definition (as we understand it):
  as $x$ travels from $A$ to $A'$
 in the complexification $V_\Ce\setminus H_\Ce$ the phase of a stable objects in $\CC_n\setminus \CC_{n+1}$ is shifted by  $n\pi$, hence the homological shift by $n$
 in requirement (2). This heuristics suggests that given a real variation of stability conditions one might expect a map from a connected covering of the complexification
 $V^0_\Ce=V_\Ce\setminus \bigcup\limits_{H\in \Sigma} H_\Ce$ to $Stab(\CC)$
 sending a point $x$ in an alcove $A$ to the stability corresponding to the central charge
 $Z(x)$ and $t$-structure $\tau_A$ (notice that the choice of $V_+$ defines a homotopy
 class of a path in $V_\Ce^0$  between any two alcoves, it is fixed by the requirement
 that the path maps $(0,1)$ to $V+ iV_+$; this defines compatible lifting of all alcoves
 to the covering).
  The main Theorem of this note is a partial result in that direction. However, the fact that
 we get a map from a covering of a proper subset in $V^0_\Ce$ which is not even
 homotopy equivalent to the whole space, and  have to use
 a somewhat unnatural quasi-exponential map is an indication of technical difficulties in connecting the two definitions. We expect even more serious difficulties in the cases when filtrations $(\CC_n)$ do not reduce to a two step filtration.

 Instead of trying to establish a direct relation between the two structures, it may be more fruitful to view them as different implementations of the same intuition of "physical" origin and possibly try to find a common generalization of the two (cf. the end of the second paragraph in \cite{Br_Se}).

\end{Rem}

\begin{Rem} In many cases (including the examples considered in this paper) 
one has natural equivalences $D^b(\AA')\cong \CC\cong D^b(\AA)$. The 
resulting equivalence $D^b(\AA')\cong D^b(\AA)$ belongs to a
class of equivalences which appeared in the work of Chuang -- Rouquier and
 Craven -- Rouquier
under the name of {\em perverse equivalences}
see \cite{ChR}, \cite{CR} and references therein
 (our setting may be slightly more general, but the generalization is straightforward).

\end{Rem}

\begin{Ex}\label{sliceex}
Let $\CC=D^b(Coh_{\BB_e}(X))$ as above, $V=\h_\RE^*$ and let $\Sigma$ consist of the affine coroot hyperplanes. Let $V^+$ be the positive Weyl chamber. Let $\tau:A\to \tau_A$ be the map described in \cite[1.8]{BM}. The polynomial map $Z:\h^*_\RE \to K^0(\CC)_\RE^*$ is characterized uniquely by its values at the points of the lattice $\Lambda\subset \h^*$; these values are given by
\begin{equation}\label{Zatla}
\<Z(\lambda), [\F]\>=\chi(\F\otimes \OO(\lambda) ) ,
\end{equation}
where $\chi$ denotes  the Euler characteristic and $\OO(\lambda)$ is the line bundle
attached to $\lambda$.

 Proposition \ref{prop1} below implies that this data provides an example
of a real variation of stability conditions, see Theorem \ref{thm2} for a stronger statement. 
Notice that in this case for every pair of neighboring
alcoves as in part 2 of the Definition the filtration $\CC_n$ is just a two term filtration, i.e. $\CC_2=\{0\}$. Another special feature of this example is that
all the $t$-structures $\tau(A)$ lie in one orbit of the group of automorphisms
of $\CC$ (in fact, of the group $B_{aff}$ acting on $\CC$). 
\end{Ex}

\begin{Rem}
\label{new_conj}
We expect 
Example \ref{sliceex} of a real variation of stability conditions 
to be a part of a richer structure. Namely, let $F\subset V$ be a face of
an arbitrary codimension and let $A$, $A'$ be two alcoves containing $F$ in 
its closure. Assume that $A$ and $A'$ are opposite with respect to $F$,
i.e. that there exists a line segment $[a,a']$ passing through $F$ with 
 endpoints satisfying
 $a\in A$, $a'\in A'$. Let $p$ be a path in the complexification
of $[a,a']$ going around $\{x\}=[a,a']\cap F$
 along a small loop in the positive 
direction. We conjecture that $\phi_p$ is a perverse equivalence governed
by the order of vanishing at $x$ of the central charge polynomials restricted
to  $[a,a']$. This can likely be restated as a real variation of stabilities
parametrized by the real points of a variety which maps birationally to $V$ so that
the preimage of $\Sigma$ is a divisor with normal crossings. 
The variety is 
closely related to
the De Concini Procesi partial compactification of 
$V^0_\Ce$ \cite{DCP}.
\end{Rem}

\begin{Rem}
Given a real variation of stability conditions one can produce a new one by adding 
an arbitrary hyperplane to the collection;
 the new set of alcoves $Alc'$ maps to the old one $Alc$, so one can define the map from 
$Alc'$ to the set of $t$-structures by composing the given map from $Alc$ with the map $Alc'\to Alc$. 
Let us say that a real variation of stability conditions is nondegenerate if it can not be obtained
this way from another one.   One can check that the real variation of stability conditions
in Example \ref{sliceex} is nondegenerate; moreover, for each two neighboring alcoves the
corresponding $t$-structures are different.
\end{Rem}

\subsection{Real variation of stabilities and automorphisms of derived
categories}\label{secrv}
In some examples in the literature
(see e.g. \cite{Br}, \cite{BrK3}, \cite{BrP2}, \cite{IUU})
(a component of) the space $Stab(\CC)$ is realized as a covering of a domain
in $K^0(\CC)_\Ce^*$ where the group of automorphisms of $\CC$ acts by deck transformations. We suggest the following counterpart of this picture in the framework of real variations of stability conditions.

\begin{Def}\label{rvarsym}
%Let us say that
 A real variation of stability conditions is {\em symmetric} if
the following holds.

\begin{enumerate} \item
For any alcoves $A$, $A'$ as in part 2 of Definition \ref{def1} there exists an auto-equivalence $m_{A,A'}$ of $\CC$ preserving the subcategories $\CC_n\subset \CC$, so that
the induced auto-equivalence of $\CC_n/\CC_{n+1}$ is isomorphic to the shift functor $M\mapsto M [2n]$.

\item The auto-equivalences $m_{A,A'}$ can be chosen so that the following holds.
 Consider the groupoid $P(V^0_\Ce)$ whose objects are alcoves and morphisms
from $A$ to $A'$ are homotopy
 classes of paths in $V^0_\Ce$ starting at $A$ and ending at $A'$.
Then there exists a functor from $F:P(V^0_\Ce)\to Cat$ where $ Cat$ is the category of categories
with morphisms being the isomorphism classes of functors,\footnote{For the sake of brevity
here 
and in Definition \ref{locsysdef} we present a weak  version
of the definition  which only deals with isomorphism classes of equivalences,
we do not address the issue of fixing the isomorphisms between the  equivalences in a
compatible way. 
One can upgrade it to a definition of a finer structure involving
the 2-category of categories, or a version of the formalism of infinity categories.} such that:

\begin{enumerate}
\item $F(A)=\CC$ for all $A\in \Alc$.

\item For $A$, $A'\in \Alc$ as above $F$
 sends the class of a path going from $A$ to $A'$ around
 $H$ in the positive direction to the identity functor.

\item For $A$, $A'$ as above $F$ sends the class
 of a path going from $A'$ to $A$ in the positive direction to $m_{A,A'}$.
\end{enumerate}
\end{enumerate}

\end{Def}

\begin{Rem}
Groupoid $P(V_\Ce^0)$ appearing above admits an explicit 
description in terms of generators and relations generalizing a standard
presentation of the (affine) braid group, see \cite{Sa}. 
\end{Rem}

Data as in Definition \ref{rvarsym}  yields the following more symmetric structure which does not involve the choice of the "positive cone" $V^+$.

\begin{Def}\label{locsysdef} Let $V,\Sigma$ be as in Definition \ref{def1}.
A {\em local system of categories with a stability condition}
 is the following collection of data. 

To an alcove $A$ one assigns a triangulated category $\CC_A$ with a bounded $t$-structure
$\tau_A$ whose heart is denoted by $\A_A$. 
%an abelian category $\A_A$, 
To every
 homotopy class
of a path $p$ connecting $A$ to $A'$ in $V^0_\Ce$ there corresponds
 a triangulated equivalence
$\phi_p: \CC_A\to \CC_{A'}$,
%D^b(\A_A)\to D^b(\A_{A'})$, 
which combine into a functor $P(V_\Ce^0)\to Cat$ (notations of Definition \ref{rvarsym}), i.e.
we have $\phi_{p p'}
\cong \phi_p\circ \phi_{p'}$ and 
$\phi_p\cong Id_{\CC_A}$ for a 
loop $p:[0,1]\to A$.

We assume that for any path $p$ connecting an alcove $A$ to itself the induced
automorphism of $K^0(\phi_p):K^0(\CC_A)\to K^0(\CC_A)$ equals identity.
We set $K=K^0(\CC_A)$ for some alcove $A$, thus $K\cong K^0(\CC_A)$ canonically
for any alcove $A$.

Finally, we assume that a polynomial {\em central charge} map
 $Z:V\to (K\otimes \RE)^*$ is fixed so that 

\begin{enumerate}
\item $\langle Z(x), [M] \rangle >0$ for $x\in A$, $M\in \A_A$, $M\ne 0$.

\item Let $A$ and $A'$ be two alcoves sharing a codimension one
face $H$. Let $p$ be the path from $A$ to $A'$ going around $H$ in the positive
direction. Then the  $t$-structure $\tau_A$ on $\CC_A$ is related
to the image of the $t$-structure $\tau_{A'}$  under
$\phi_p^{-1}$ as spelled out in
%$\phi_p$ is a perverse equivalence characterized in terms
%of the central charge map as in 
Definition \ref{def1}(2). 
\end{enumerate}
\end{Def}
 
 As was pointed out above, a symmetric real variation of stabilities 
 yields a local system of categories with a stability condition by setting $\CC_A:=\CC$,
 $\tau_A:=\tau(A)$ etc.
 
Conversely, it is easy to see that given a  local system of categories with a stability
condition together with an additional choice of a cone
$V^+$ as in Definition \ref{def1}
 one can obtain a symmetric real variation of stabilities in the following
way. For any two alcoves $A,A'$ there exists a unique homotopy class of a path
$p_{A,A'}:[0,1]\to V^0_\Ce$ such that $p:(0,1)\to V+\sqrt{-1}V_+$; moreover,
we have $p_{A,A''}\sim p_{A',A''} p_{A,A'}$ for any triple of 
alcoves $A,\,A',\, A''$. Thus we can fix a triangulated category $\CC$
together with equivalences $\phi_A:\CC\to D^b(\A_A)$, so that
$\phi_{A'}\phi_A^{-1}\cong \phi_{p_{A,A'}}$. The category $\CC$ then comes
equipped with the data described in the definition of a symmetric
 real variation of stabilities.

%The auto-equivalence $m_{A,A'}$ can be thought of as monodromy along the loop in $V^0_\Ce$ starting and ending at $A$ and running around the common face of $A$ and $A'$ in the positive direction.

\subsection{Symmetric real variation of stabilities and symplectic resolutions}
\label{secsymp}
Let $\pi: X\to Y$ be a conical (homogeneous) symplectic resolution over the field $k$. Recall that this means
that $\pi$ is a resolution of singularities, $X$ carries an algebraic symplectic form $\omega$
and the multiplicative group $\Gm$ acts on $X$, $Y$ compatibly, dilating the symplectic form
so that $t^*(w)=t^2\omega$. The action of $\Gm$ on $Y$ is assumed to contract $Y$ to a point
$0\in Y$; this implies that $Y$ is affine. We set $\CC=D^b(Coh_{\pi^{-1}(0)}(X))$.

Along with $\CC$ we will also need to consider its graded version
$\CC_{gr}:=D^b(Coh^{\Gm}_{\pi^{-1}(0)}(X))$. Notice that $K^0(\CC_{gr})$ is module
over $\Zet[q,q^{-1}]$ where $q$ acts by twisting with the tautological character of $\Gm$.  
For $t\in \Ce^\times$ we set $K(\CC_{gr})_t=K(\CC_{gr})\otimes_{\Zet[q,q^{-1}] }\Ce$, where $\Zet[q,q^{-1}]$ maps to $\Ce$ by sending $q$ to $t$.

Assume for simplicity that $Pic(X)\cong \Zet^r$. 
Set $V=Pic(X)\otimes \RE$; let $V^+$ is spanned by the ample cone, and let central charge map $Z$ be given by  formula \eqref{Zatla}.

We now state a conjecture due to the second author and A.~Okounkov, see also discussion
in \cite[\S 1.10]{BMO}, \cite[\S 1.1.6]{MO} etc. 

In the statement we need  some basic information about quantizations in positive characteristic.
It is known that for $k$ of characteristic zero, hence also of large positive
characteristic, $H^i(X,\O)=0$ for $i>0$, thus $X$ is admissible in the sense of 
\cite[Definition 1.21]{BK1}. Thus the main result of \cite{BK1} implies existence of a canonical quantization of $X$ over $k$ of  characteristic $p>0$. Furthermore,  \cite[Theorem 1.23]{BK1}
and discussion following its proof show that for $\lambda\in Pic(X)$ one can twist
the canonical quantization by $\lambda$, we call the result the quantization with parameter 
$\lambda$. All these quantizations are Frobenius constant (see \cite{BK}, \cite{BK1}
for the definition) thus the result is an Azumaya algebra on the Frobenius twist $X^{(1)}$ of $X$.
We denote this Azumaya algebra by $\A_\la$, we also fix an isomorphism of schemes $X\cong X^{(1)}$. Notice that the isomorphism class of $\A_\la$ depends only on the image of $\lambda$
in $Pic(X)/pPic(X)$; however, we view $\lambda$, not that image, as a parameter for quantization,
this allows us to fix a compatible system of Morita equivalences $A_\la\sim A_\mu$, 
$\la,\mu\in Pic(X)$, these Morita equivalences are used in making sense of compatibility
for splitting of Azumaya algebras in part (3) of the following Conjecture.

\begin{Conj}\label{conj}
The above data of $\CC$, $V$, $Z$, $V^+$ admits a natural extension to a symmetric 
real variation of stability conditions  
with the following properties.

\begin{enumerate}
\item The action of $\pi_1(V_\Ce^0)$ on $\CC$ lifts to an action on 
$\CC_{gr}$.

\item Assume that $k=\Ce$. Set $T=N(X)\otimes \Ce^*$ and set 
$$D=\bigcup\limits_{H\in \Sigma} \{ \exp(2\pi i h)\ |\ h\in H\}.$$ 
Then the (small) {\em $\Ce^*$-equivariant quantum cohomology of $X$}
defines a  family of flat connections on the trivial vector bundle
over $T$ with fiber $H^*(X)$ depending on a parameter $h$ (the equivariant parameter).
These connections have regular singularities on $D$.

We have $\pi_1(V_\Ce^0)\subset \pi_1(T\setminus D)$. For a given generic
value $h$ of the parameter  the monodromy of the corresponding connection 
restricted to $\pi_1(V_\Ce^0)$ is isomorphic to the induced action on $K(\CC_{gr})_{\exp(2\pi ih)}$.

\item Assume that $k$ is a field of characteristic $p\gg \dim(X)$. 
Suppose that $\lambda\in Pic(X)$ is such that $\frac{\lambda}{p}\in A$.

The Azumaya algebra $\A_\la$ is split on the formal neighborhood of $\pi^{-1}(0)$,
we fix such splitting for all $\la\in Pic(X)$ in a compatible way.  Thus we get an equivalence between
$\CC$
 and the category $\A_\la-mod_{\pi^{-1}}(0)$ of coherent sheaves of  $\A_\la$-modules
supported on $\pi^{-1}(0)$. The composed functor
$$\CC\cong D^b(\A_\la-mod_{\pi^{-1}(0)} )\overset{R\Gamma}{\longrightarrow} D^b(\Gamma(\A_\la)-mod)$$
is an equivalence of triangulated categories sending $\tau_A$ to the tautological
$t$-structure. 

\item Let  $X_R\to Y_R$ be a conical symplectic resolution  
over a commutative ring $R$ finitely generated over $\Zet$.
Then  possibly after a finite localization of $R$,   the symmetric
real variations of stabilities from (2,3) exist for all
base changes $X_k\to Y_k$ for an algebraically
closed field $k$, and data for different $k$ are
   compatible as follows.

The set $\Sigma$  does not depend on $k$, 
and for each alcove $A$ 
there exists
%a symmetric real variation of stabilities
a $t$-structure $\tau_A^R$
on $D^b(Coh_{\pi^{-1}(0)}(X_R))$ compatible under base change with the $t$-structure
$\tau_A$ on $D^b(Coh_{\pi^{-1}(0)}(X_k))$.
\end{enumerate}
\end{Conj}

\begin{Rem}
The idea that
stability conditions should be related to quantum cohomology goes back to Bridgeland (see e.g.
\cite[\S 7.2]{Br_Se})
who proposed it based on heuristics of mirror symmetry. However,
our context differs from that considered by Bridgeland in several  aspects;
in particular,  it is essential for us to work with {\em equivariant} quantum cohomology, while 
there seems to be no available description  of the (conjectural) mirror dual
counterpart of equivariant quantum cohomology of a local Calabi-Yau; here
we refer to equivariance with respect to an action of the multiplicative group dilating the volume form. 
 \end{Rem}

\begin{Thm}\label{thm2} Let $X\to Y$ be as in section \ref{sec1} and assume that
the map $H^2(G/B)\to H^2(X)$ is an isomorphism.\footnote{This condition
holds quite often, in particular it always holds if (the Dynkin graph of) $G$ is simply-laced;
see \cite[Theorem 1.3]{list} for the list of exceptional cases.}
Then Conjecture \ref{conj} holds.
\end{Thm}

\proof Existence of a symmetric real variation of stability conditions with $\Sigma$
being the collection of affine coroot hyperplanes, as well as its relation to quantization in positive
characteristic described in property (3), follow from Proposition \ref{prop1} and its proof.  
The lifting of the $B_{aff}$ action to the equivariant category is addressed in \cite{BM}, this
yields part (1).
The $t$-structures
are constructed in \cite{BM} for slices to nilpotent elements defined over $\Zet[\frac{1}{h!}]$
where $h$ is the Coxeter number, so part (4) is also clear. Finally, 
equivariant quantum cohomology of $X$ was computed under the above assumptions in 
\cite{BMO}. The monodromy of the resulting connection (called the affine KZ, or trigonometric Dunkl)
connection is isomorphic to the $B_{aff}$ module coming from the standard module
for the affine Hecke algebra by a result of \cite{Che}, see  \cite[Proposition 2.3]{BMO}. This is isomorphic to $K^0(\CC_{gr})$ by \cite[Theorem 1.3.2(b)]{BM}. \epf

\section{Proof of Theorem \ref{thm1}}
\subsection{Positivity property}

\begin{Prop}\label{prop1}
 Let $A\in Alc$ and let $M\ne 0$ be an object in the heart of $\tau_A$. Let $A'\in Alc$ be a neighboring alcove separated from $A$ by a codimension one face $F$.

a)
 The function $d_M:x\mapsto \langle E(x), [M] \rangle$ is a polynomial taking positive real values on $x\in A$.

b) Either $d_M$ takes positive values on $F$, or $d_M$ has a zero of order one on $F$.

c)  Assume that  $M$ is irreducible.
If  $d_M|_F =0$, then $b_{A,A'}(M)=M[\pm 1]$, where  the $+$ (respectively, $-$) sign should be taken if $A$ lies above (respectively, below) $A'$. 

Otherwise $b_{A,A'}(M)=\tilde M$, where $\tilde M$
is an object in the heart of $\tau_A$
which fits into one of the two exact sequences:
$$0\to M' \to \tilde M\to M\to 0,$$
$$0\to M\to \tilde M\to M'\to 0,$$
 where the first (respectively, the second) exact sequence applies if $A$ lies above (respectively,
 below) $A'$. Furthermore, $d_{M'}|_F=0$ and $M$ is the only simple quotient (respectively, sub)
 module of $\tilde M$.
\end{Prop}

{\em Proof.} In this proof we will work over an algebraically closed field $k$ possibly of characteristic
$p>0$.
For large $p$  conjugacy classes of
nilpotent elements in $\g_k$ are in a natural bijection with those in $\g=\g_\Ce$, 
we fix $e_k\in \g_k$ in the class corresponding to the class of $e$. 

The construction
 of  $t$-structures $\tau_A^k$, $A\in Alc$ on $\CC_k=D^b(Coh_{\BB_{e_k}}(X_k))$ (in self-explanatory notation) is carried out in \cite{BM} for all $k$ except those of positive characteristic $p \leq h$, where $h$ is the Coxeter number of $G$.
 
 Let $\A_A^k$ denote the heart of $\tau_A^k$. Then for $p\gg 0$ the set of irreducible
 object in $\A_K^k$ is in a natural bijection with irreducible objects of $\A_A^\Ce$ 
 (see \cite[\S 5.1.4]{BM}), so that
 the classes of corresponding objects match under the standard
 identification   $K^0(\BB_{e_\Ce})=K^0(\BB_{e_k})$ explained e.g. in
 \cite[\S 7]{BMR}.

Fix  $p\gg 0$, and  let $\la\in \La$ be such that $\frac{\la+\rho}{p}\in A$,
where $\rho$ is the sum of fundamental weights. Then it is shown in \cite[\S 6]{BMR}
that 
\begin{equation}\label{dim}
p^{\dim \BB} d_\F(\frac{\la+\rho}{p})= \dim (\Gamma_\la(\F)),
\end{equation}
where  $\F\in \CC_k$ and  $\Gamma_\la:\CC_k\to D^b(\MM_\lambda)$ is 
%an equivalence of triangulated categories
a functor sending $\A_A^k$ to $\MM_\la$,  inducing a full embedding  $\A_A^k\imbed  \MM_\la$.
 Here $\MM_\lambda$
is the category of modules over the Lie algebra $\g_k$ with 
the fixed generalized central character $(\la, e)$ (here $\la$ and $e$ are, respectively,  
the  Harish-Chandra and the Frobenius central characters, 
%and the action of  Frobenius central character factoring
%through the ring of functions on the slice $S_k$, with a power of the maximal ideal corresponding to $e$ acting by zero  (
see \cite[\S 1]{BMR} for terminology). 

Thus for a large prime number $p$ and $M\in \A_A^\Ce$, the polynomial $d_M$ takes 
a positive value $p^{-\dim \BB} \dim(\Gamma_{p\la-\rho}(M_k))$  at points $\la\in A$ such that
$p\la\in \La$; here $M_k\in \A_A^k$ is an object  whose class in the Grothendieck group matches that of $M$. 
Since the set of such points, for varying $p$, is dense, we see that $d_M(\la)\geq 0$ for 
any $\la\in A$. It remains to see that the inequality is strict.

Suppose that $d_M(\la_0)=0$ for some $\la_0\in A$. We claim that $d_M$
 is a harmonic polynomial, i.e.
 $ d_M(\la)=\frac{1}{|W|} \sum_w d_M (\la+ w(\mu))$ for all $\la$, $\mu\in \h^*$.
  This follows from the
the fact that $d_M$ has the form 
%\begin{equation}\label{form}
$d_M(\lambda) = \langle \xi, \exp(\lambda) \rangle,$
%\end{equation}
where $\xi$ is a linear functional on $Sym(\h^*)/ (Sym(\h^*)^W_+) \cong H^*(G/B)$
and $\exp:\h^*\to Sym(\h^*)/ (Sym(\h^*)^W_+)$ is the exponential map.
It is clear  that a 
harmonic polynomial which vanishes at a point but
takes non-negative values on some neighborhood of that point is identically zero. This proves (a).

The proof of (b) is based on "singular localization" Theorem of \cite{BMR2}. Namely, let $\la\in \La$ be such that $\frac{\la+\rho}{p}$ lies on the boundary
% a codimension one face 
of an alcove $A$.
 Then we still   have %$p^{\dim \BB} d_\F(\frac{\la+\rho}{p})=\dim \Ga_\la(\F)$
  the functor 
 $\Gamma_\la:\CC_k\to D^b(\MM_\la)$ which sends $\A_A^k$ into $\MM_\la$ but now it is
  not necessarily conservative (i.e. it may kill some non-zero objects). The equality \eqref{dim}
  still holds.
  
  Furthermore, the set of irreducible objects $\A_A^k$ killed by $\Gamma_\la$ depends
  only on the face $F$ containing $\frac{\la+\rho}{p}$; it is also independent of $p\gg 0$,
  where we use the above identification of the set of irreducible objects in $\A_A^k$ for varying $k$.
  %We denote 
  
  %For $\la_1,\,\la_2$ such that  $\frac{\la_i+\rho}{p}\in F$ we again
 %have an equivalence given by the translation functor $\Gamma_{\la_1,\la_2}$ satisfying
 %\eqref{tra}. 
 %
% For $\lambda_1,\, \lambda_2\in \Lambda$ with  $\frac{\la_i+\rho}{p}\in A$ we have a canonical equivalence
%({\em not} preserving dimensions of modules) called the {\em translation functor} \cite[\S 2.2]{BMR2}
%$T_{\lambda_1,\lambda_2}: \MM_{\lambda_1}\cong \MM_{\lambda_2}$
%such that
% \begin{equation}\label{tra}
%\Gamma_{\lambda_2}\cong T_{\lambda_1,\lambda_2}\circ \Gamma_{\lambda_1}.
%\end{equation}

 Assuming that $p$ is large enough, we see that if $\Gamma_\la(M)=0$ for some $\frac{\la+\rho}{p}\in F$,
then the polynomial $d_M$ vanishes at all points of $F$ with sufficiently  large prime
denominator, hence $d_M|_F\equiv 0$.

Otherwise $d_M$ takes positive values at all points of $F$ with a large prime denominator, then an argument involving harmonic polynomials as in the proof of part (a) shows it takes positive values at all points of $F$. 

It remains to see that the order of vanishing of polynomial $d_M$ on $F$ can not be greater
than one.  We claim that this is true for any harmonic polynomial: if a harmonic polynomial
$P$ has zero of order two a more on a coroot hyperplane $F$, we can apply a differential
operator with constant coefficients to $P$ to get a harmonic polynomial  which has the
form $\alpha^2Q_0$, where $\alpha$ is the equation of the hyperplane and $Q_0|_F\ne 0$.
If $C$ is the Laplace operator (the $W$-invariant order two operator with no constant term),
then $C(\alpha^2Q_0)|_F=C(\alpha^2)Q_0|_F\ne 0$, which contradicts the fact that $W$-invariant 
differential operators with constant coefficients and no constant term annihilate harmonic polynomials.

c)
Consider first  the situation over $k$ of large positive characteristic. 
In this case  the braid group action admits a convenient description in terms of
 translation functors which we presently recall.

%We have the following facts.

\begin{Claim}\label{cl}
Fix $\la,\,\mu \in \La$ with $\frac{\la+\rho}{p}\in A$, $\frac{\mu+\rho}{p}\in \bar{A}\setminus A$, where
$\bar{A}$ is the closure of $A$. 

a) $B_{aff}$ acts on $D^b(\MM_\la)$, so that $\Gamma_\la$ is compatible with the braid group 
action.

b) We have a bi-adjoint pair of exact functors called {\em translation functors}

$T_{\la\to \mu}:\MM_\la\to \MM_\mu$, $T_{\mu\to \la}:\MM_\mu\to \MM_\la$.

%The composition of adjunction arrows $L\to T_{\mu\to \la}T_{\la\to\mu}(L)\to L$
%is zero for any irreducible object $L\in \MM_\la$.
The functor $T_{\la\to\mu}$ sends an irreducible object to an irreducible one or zero,
the nonzero images of nonisomorphic irreducibles  are not isomorphic.

c) We have $\Gamma_\mu\cong T_{\la\to\mu}\circ \Gamma_\la$.

d) We now assume that $\frac{\mu+\rho}{p}$ lies in the codimension one face
$F$ separating $A$ from $A'$.
 Set $R=T_{\mu\to \la}T_{\la\to\mu}$.

Then we have  functorial exact triangles
\begin{equation}\label{down}
Id\to R\to b_{A,A'}\to Id[1]
\end{equation}
if $A$ is above $A'$;
\begin{equation}\label{up}
 b_{A,A'}\to R\to Id \to b_{A,A'}[1]
\end{equation}
if $A$ is below $A'$; here the map to/from $Id$ from/to $R$ is the adjunction arrow.

\end{Claim}

{\em Proof} \  of the Claim. The $B_{aff}$ action is introduced in   
\cite[Theorem 2.1.4, Corollary 2.1.6 ]{BMR2}, compatibility with the $B_{aff}$ action is \cite[Proposition 1.6.4]{BM}, for proof see \cite[\S 5]{Ri}: this proves (a). 

Translation functors and their adjointness are discussed  in this setting in \cite[\S 2.2.1]{BMR2}.
The properties of translation of an irreducible module to the wall are stated in \cite[\S 2.2.6, Remark 1]{BMR2}, the proof is similar to that of the corresponding fact in characteristic zero,
see \cite[\S 2.5]{BeG}, this yields (b).

Part (c) is \cite[Lemma 2.2.3(a)]{BMR2}.

It is easy to see that in notation of \cite{BM} we have $b_{A,A'}=\tilde s_\alpha^{\pm 1}$
if $A$ and $A'$ are separated by a face of type $\alpha$; here the plus sign is chosen
if and only if $A$ is above $A'$.
Thus part (d) follows by comparing \cite[Lemma 2.2.3(c)]{BMR2} and characterization of the affine braid
group action in \cite[Theorem 2.1.4, Corollary 2.1.6 ]{BMR2}. \epf

\medskip

Now we are ready to prove part (c) of the Proposition. We assume that $A$ is above $A'$, the other case is treated similarly. Assume $L\in \A_A^\Ce$ is such that $d_L|_F=0$. Let  $L_k$ be the corresponding
irreducible in $\A_A^k$. We have $\Gamma_\mu(L_k)=0$, thus, setting $M_k=\Gamma_\la(L_k)$ we get 
$R(M_k)=0$
and \eqref{down} shows that $b_{A,A'}(M_k)\cong M_k[1]$.
Since $\Gamma_\la$ is fully faithful, we see that $b_{A,A'}(L_k)=L_k[1]$.

Now assume that $d_L|_F\ne 0$, the proof of part (b) of the Proposition shows that $T_{\la\to\mu}(M_k)\ne 0$. 
Since $M_k$ is irreducible and the adjunction map $M_k\to R(M_k)$ is nonzero, 
the map $M_k\to R(M_k)$ in  \eqref{down} is injective, thus $b_{A,A'}(M_k)\in \MM_\la$ and we have
a short exact sequence in $\MM_\la$:
$$0\to M_k\to R(M_k)\to b_{A,A'}(M_k)\to 0.$$

%Applying the functor $T_{\la\to \mu}$ to this sequence we get, in view of
 %5) a short exact sequence 
 %$$0\to T_{\la\to\mu}(M_k)\to  T_{\la\to \mu}(R(M_k))\to T_{\la\to\mu}(M_k)\to 0.$$

For any irreducible $N\in \MM_\la$, $N\not\cong M_k$ we have 
\begin{equation}\label{Homo}
Hom(b_{A,A'}(M_k),N)\subset
Hom(R(M_k),N)=Hom(\Gamma_{\la\to\mu}(M_k),\Gamma_{\la\to\mu}(N))=0
\end{equation}
  in view of Claim \ref{cl}(b). Thus there exists a nonzero map $f:b_{A,A'}(M_k)\to M_k$.
  Applying $\Gamma_{\la\to\mu}$ to $f$ we get a nonzero, hence surjective
  map from $b_{A,A'}(M_k)$ to the irreducible (by Claim \ref{cl}(b)) module $\Gamma_\la(M_k)$.
  Now Claim \ref{cl}(c) implies that $\Gamma_{\la\to\mu}(f)$ is an isomorphism, hence
  $\Gamma_{\la\to\mu}(Ker(f))=0$.   
  
Taking into account \eqref{Homo}
we get that $M_k$ is the only irreducible quotient of $b_{A,A'}(M_k)$, while
 the kernel of the map $b_{A,A'}(M_k)\to M_k$  is killed by $T_{\la\to \mu}$. 
 Since $\Gamma_\la$ is fully faithful, we see that
  $L_k$ is the only irreducible quotient of $b_{A,A'}(L_k)$,
 while the kernel $K$ of the map $b_{A,A'}(L_k)\to L_k$ satisfies  $d_K|_F=0$.
It follows that the same is true for $L$. 
\epf

\subsection{End of proof of the Theorem.} \label{sect3}
We start with an auxiliary Lemma, probably well known to the experts.

To state it we recall notations of \cite{Br1}. Given a stability condition $\sigma$ on a triangulated
category $\CC$ one has full subcategories $\PP(\phi)$ and $\PP(I)$ in $\CC$ for $\phi\in \RE$
and an interval $I\subset \RE$ (we will also write $\PP_\sigma$ when we need to emphasize dependence on $\sigma$). The categories $\PP((a,a+1])$ are abelian (in fact, they are
hearts of bounded $t$-structures on $\CC$), while the categories $\PP((a,b))$ for $a<b<a+1$
are quasi-abelian. The categories $\PP(\phi)$ are also abelian \cite[Lemma 5.2]{Br1} and
 simple objects of these abelian categories are called stable.

\begin{Lem}\label{LBr} a) Given a locally finite stability condition on a triangulated
category $\CC$ and an object $M\in \CC$ the following are equivalent.

\begin{enumerate}

%\item $M$ is a simple object in the abelian category $\PP((a,a+1])$ for some $a \in \RE$.

\item $M$ is a simple object in the quasi-abelian category $\PP((a,b))$ for some $a$, $b$, 
$a<b<a+1$.

\item $M$ is stable.

\end{enumerate}

b) Suppose that $M$ is stable  with
respect to a locally finite stability condition $\sigma$. Then there exists an open neighborhood
$U$ of $\sigma$ in $Stab(\CC)$ such that $M$ is stable with respect to any $\sigma'\in U$.
\end{Lem}

\proof a) 
%The implication $(1) \Rightarrow (2)$ is clear from $\PP((a,a+1])=\bigcup\limits
%_{b\in (a,a+1)} \PP(b,a+1)$ which shows that an object 
%
 It is clear from definitions that (1) implies (2). To check that (2) implies (1)
assume that  $M$ is stable of phase $t$. By the definition of local finiteness 
we can find a finite length quasi-abelian category $\BB = \PP((t-a,t+a))$ ($a<\frac{1}{2}$).
 The object $M$ has a finite Jordan-Hoelder series in $\BB$.
  Thus there are only finitely many elements in $K^0(\CC)$ which
can be represented by a subobject of $M$ in $\BB$, since
such a class is a sum of classes of some of the simple constituents.
Since $M$ is stable, the class of  a subobject $N\ne M$ has phase strictly less than $t$. Thus there exists  $s\in (t-a,t)$ such that the class of any such subobject
has phase less than $s$. Then $M$ is irreducible 
in $\BB'=\PP((s, t+a))$. 

b)  Using (a) we find $a$, $b$, $a<b<a+1$, such that $M$ is simple in $\PP_\sigma((a,b))$. 
Fix $\alpha$, $\beta$, so that $a<\alpha< \phi <\beta<b$, where $\phi$ is the phase of $M$. By the definition of topology
on the set of %locally finite 
stability conditions  (cf. \cite[Lemma 6.1]{Br1}), there exists an open neighborhood
$U$ of $\sigma$ such that 
$M\in \PP_{\sigma'}((\alpha,\beta))\subset \PP_\sigma((a,b))$ for $\sigma'\in U$.
Thus $M$ is a simple object in $\PP_{\sigma'}((\alpha, \beta))$ for $\sigma'\in U$, applying
part (a) again we see that $M$ is stable with respect to such $\sigma'$. 
 \epf
 
\bigskip

We now complete the proof of the Theorem.
We  construct the required map as follows. Recall that $A_0$ denotes the fundamental alcove. It is easy to see that the set
$$S=\{(\la,\mu) \in \h^*_\RE\times \h^*_\RE\cong \h^* \ | \ (\la\in A_0) \bigvee (\la \in \bar{A_0}, \, \mu \in A_0) \}$$
is a fundamental domain for the action of $W_{aff}$ on $V^{reg}$, here we used identification \eqref{xiy}.
 In fact, $S$ is the intersection of a contractible fundamental domain for the action of $W_{aff}$ on
  $\h^*$ with $V$.

Thus  a point in $\wt{V^{reg}}$ can be represented by a pair $(b,x)$ where $x\in S$ and $b$ is a homotopy class of path
from $A_0$ to some alcove $A\in Alc$ (the projection to $V^{reg}$ is then given by $(b,x)\mapsto \bar{b}(x)$ where $\bar{b}$
is the element of $W_{aff}$ corresponding to $b$). We define the map $\iota$ by:
$$\iota: (b,x)\mapsto {\mathfrak S}(b(\tau_{A_0}), \ii E(\bar{b}(x))),$$
where we use the same notation for $b$ and the corresponding element of $B_{aff}$. 
Here $\mathfrak S$ denotes the stability condition obtained from the given $t$-structure
and central charge by means of  \cite[Proposition 5.3]{Br1}.
 Recall from section \ref{11} that the heart of $\tau_A$ is a finite length abelian category,
it is easy to deduce that any $\sigma\in \iota(\wt{V^{reg}})$ is locally finite in the sense
of \cite[Definition 5.7]{Br1}, thus the map lands in the space $Stab(\CC)$ of locally finite stability conditions.

It is clear from the definition that the map $\iota$ is $B_{aff}$ equivariant.
It remains to check that $\iota$ is continuous.  
% This is done using Lemma \ref{LBr}(b).

 We have to check continuity at the boundary of the region corresponding to a given $b\in B_{aff}$. Without loss of generality we can assume $b=1$. Let $\bla\in S$ be a boundary  point
 (i.e. $\bla$ is not an inner point  of $S$). Thus $\bla=(\la,\mu)$, $\la\in F$, $\mu\in A_0$,
 where $F$ is a codimension one face of $A_0$. The face $F$ belongs to a unique hyperplane
which corresponds to an affine coroot $\al$ attached to a vertex of the affine Dynkin diagram.
 It is easy to see that a small neighborhood of $(1,\bla)$ is contained in the union of regions
 corresponding to $1$ and $\tilde s_\al^{-1}$ (notice that $b_{A,A'}=\tilde s_\al^{-1}$ if $\al$ corresponds
 to a vertex of the finite Dynkin graph, in which case $A'$ is below $A$; and $b_{A,A'}=
 \tilde s_\al$ otherwise\footnote{Here we use a (possibly nonstandard) identification between
 $\pi_1((\h^*)^{reg}/W_{aff})$ and the Artin group given by the standard presentation: we send
 the generator $\tilde s_\alpha$ to the loop going around the coroot hyperplane in the {\em negative}
 (clockwise) direction. This allows to avoid extra notational complication in the statement of Theorem
 \ref{thm1}, while using the $B_{aff}$ action of \cite{BR}. }
% This choice is also related to the choice made by Bridgeland in his definition of 
 %the $t$-structure corresponding to a given stability condition: he requires the phases of the 
 %objects in the heart to lie in $(0,1]$; one could alternatively make it lie $[0,1)$, then
 %Theorem \ref{thm1} would hold with the other choice of generators for the fundamental group.}).
   
 By \cite[Theorem 1.2]{Br1}  there exists a neighborhood $\tilde U$ of the point $(\tau_{A_0},\bla)$ in $Stab(\CC)$ mapping isomorphically to a neighborhood $U$ of $\ii E(\bla)$ in
 $K^0(\CC)^*$. It suffices to see that for a small enough $U$ and a point $\tilde z
 \in \tilde U$  mapping to $z\in s_\al(S)\cap U$, the $t$-structure underlying $\tilde z$ is $\tau_{A'}=  \tilde s_\al ^{-1}(\tau_A)$; here $A'$ is the neighboring alcove, $A'=s_\al(A_0)$.

 Let $M$ be an irreducible object. It suffices to show that $\tilde s_\al^{-1}
 (M)$ lies in the heart of the $t$-structure of $\tilde z$. The proof is similar to the proof of \cite[Lemma 3.5]{Br}.

 Consider the dichotomy of Proposition \ref{prop1}(c). 
 If $d_M|_F=0$, then according to Proposition \ref{prop1}(c) $%b_{A,A'}
\tilde s_\al ^{-1}(M)=M[-1]$. Using  Lemma \ref{LBr}(b) we see that $M$ is stable for $\tilde z$ (if $\tilde U$ is small enough).
We have $\langle Z(\bla),[M]\rangle\in \RE_{<0}$, while $\langle z, [M] \rangle$ lies in the open 
 lower halfplane: the latter fact follows from the order one vanishing statement in Proposition 
 \ref{prop1}(b). Thus the phase of $M$ with respect to $\tilde z$ is in $(1,2)$,
 so  $M[ -1]$ is in the heart of the $t$-structure since its phase with respect to $\tilde z$
is in $(0,1)$.
 
  Or else $\tilde s_\al^{-1}(M)=\tilde M$ lies in the heart of $\tau_A$,
%its Harder-Narasimhan filtration has length two, and both stable subquotients remain stable and have phases in $[0,1)$ in stability $\tilde z$. 
%
%and all its stable subquotients with respect to $\tilde z$ have 
it is stable in stability $\tilde z$ and has phase in $(0,1)$. To see
 this recall that $\langle Z(\bla), [M]\rangle$ lies in the open upper halfplane,
 while $\langle Z(\bla), [L] \rangle \in \RE_{<0}$; since every nonzero proper subobject 
 $N$ in $ \tilde M$ contains $M$ we see that
 $$\langle Z(\bla), [N] \rangle = \langle Z(\bla), [\tilde M ] \rangle + s$$
   for some $s\in \RE_{>0}$. Thus
   the phase of $N$ is smaller than the phase of $\tilde M$, i.e. $\tilde M$ is stable.
    It follows that the same is true for stability conditions
in a neighborhood of $\tilde z$. 
\epf

 \end{document}